\documentclass[a4paper,12pt]{article}
\usepackage{amssymb}
\usepackage{amsthm}
\theoremstyle{plain}

\newtheorem{Proposition}{Proposition}
\newtheorem{Theorem}{Theorem}

\theoremstyle{definition}

\newtheorem{Remark}{Remark}
\newcommand{\reg}{\mathrm{reg}}
\newcommand{\im}{\mathrm{Im}}
\newcommand{\Ker}{\mathrm{Ker}}
\newcommand{\codim}{\mathrm{codim}}
\newcommand{\PP}{\mathbb{P}}
\newcommand{\C}{\mathbb{C}}
\newcommand{\K}{\mathbb{K}}
\newcommand{\N}{{\mathbb N}}
\newcommand{\R}{{\mathbb R}}

\begin{document}
\title{On the osculatory behaviour of \\ 
higher dimensional projective varieties}
\author{Edoardo Ballico and Claudio Fontanari}
\date{}
\maketitle 

\begin{small}
\begin{center}
\textbf{Abstract}
\end{center}
Here we explore the geometry of the osculating spaces to
projective varieties of arbitrary dimension. 
In particular, we classify varieties having very degenerate 
higher order osculating spaces and we determine mild conditions
for the existence of inflectionary points. 
\vspace{0.5cm}

\noindent AMS Subject Classification: 14N05.

\noindent Keywords: bundle of principal parts, osculating spaces.
\end{small}
\vspace{0.5cm}

\section{Introduction}
The geometry of osculating spaces to projective varieties is a very 
classical but still widely open subject. As it is well-known, while 
the dimension of the tangent space at a smooth point is always equal 
to the dimension of the variety, higher order osculating spaces can 
be strictly smaller than expected also at a general point. 

The investigation of algebraic surfaces having defective second order 
osculating space was inaugurated in 1907 by Corrado Segre in his 
seminal paper \cite{Segre:07}. The early developments in the field 
are witnessed by a long series of nice contributions, among which 
we cannot resist to mention at least \cite{Terracini:12} by Alessandro 
Terracini, \cite{Bompiani:19} by Enrico Bompiani, \cite{Togliatti:29} 
and \cite{Togliatti:46} by Eugenio To\-gliat\-ti. We point out that the 
work of Togliatti has been recently reconsidered by several authors: 
the main result of \cite{Togliatti:46} is reproved in \cite{ilardi} 
and \cite{fraila}, while the crucial example in \cite{Togliatti:29} 
is analized in \cite{perkinson} and \cite{lanmal}. However, the modern 
approach to the subject is mainly concerned with the osculatory behaviour 
at every (not just at the general) point: see for instance \cite{shifrin}, 
\cite{fkpt}, \cite{pt}, \cite{bpt}, \cite{lan}. 

Here we address the case of projective varieties of arbitrary dimension 
from various points of view. Namely, in Section~\ref{s1} we apply the 
classical approach of Bompiani and we obtain a rough classification of 
projective varieties having very degenerate higher order osculating 
spaces (see Theorem~\ref{main}). In Sections~\ref{s2} and \ref{s3}, 
instead, we apply the modern theory of vector bundles and we find 
explicit conditions assuring the existence of inflectionary points on 
projective varieties (see Theorems~\ref{a1} and \ref{a5}). 
Other related results are collected in Propositions~\ref{a3}, \ref{a4}, 
and \ref{lanteri}.  

This research is part of the T.A.S.C.A. project of I.N.d.A.M., supported by 
P.A.T. (Trento) and M.I.U.R. (Italy).

\section{}\label{s1}
Let $X \subset \PP^r$ be an integral nondegenerate projective 
variety of dimension $n$ defined over the field $\C$. 
Let $p \in X$ and consider a lifting 
\begin{eqnarray*}
U \subseteq \C^n &\longrightarrow& \C^{r+1} \setminus \{ 0 \} \\
t &\longmapsto& p(t)
\end{eqnarray*}
of a local regular parametrization of $X$ centered in $p$. 
The osculating space $T(m,p,X)$ of order $m$ at $p \in X$ 
is by definition the linear span of the points $[p_I(0)] \in \PP^r$, 
where $I$ is a multi-index such that $\vert I \vert \le m$. 
We say that $X$ satisfies a differential equation of order $m$ at $p$ if 
$$
\sum_{\vert I \vert \le m} a_I(t_1, \ldots, t_n) p_I(t_1, \ldots, t_n) = 0
$$ 
in $U$, with $a_I(t_1, \ldots, t_n) \ne 0$ for some $I$ with 
$\vert I \vert = m$. Hence we have 
$$
\dim T(m,p,X) = {{n+m}\choose{n}} - N - 1
$$
where $N$ is the number of independent differential equations of order $m$ 
satisfied by $X$ at $p$. In \cite{bf} we proved the following generalization 
of a classical result by Bompiani: 

\begin{Theorem} \label{bompiani} \emph{(\cite{bf})}
Let $X \subset \PP^r$ be a smooth variety and let $p \in X$ be a general 
point. Assume that $\dim T(m,p,X) = h$ and $\dim T(m+1,p,X) = h + k$ with 
$1 \le k \le n-1$. 
Then either $X \subset \PP^{h+k}$ or $X$ is covered by infinitely 
many subvarieties $Y$ of dimension at least $n-k$ such that 
$Y \subset \PP^{h-m}$.
\end{Theorem}

As an application in the spirit of Bompiani's paper \cite{Bompiani:19}, 
we are going to classify smooth projective varieties satisfying many 
differential equations at a general point.

\begin{Remark}\label{hopf}
Fix $p \in X$ and let $V_m(p)$ be the $\C$-vector space of the differential 
equations of order $\le m$ satisfied by $X$ at $p$. Set $W_m(p) = 
V_m(p) / V_{m-1}(p)$ and define $h_m(p) := \dim W_m(p)$. There is a natural 
linear map 
\begin{eqnarray*}
H^0(\PP^{n-1}, \mathcal{O}(d)) \otimes W_m(p) &\longrightarrow& W_{m+d}(p) \\
\left( \frac{\partial}{\partial x_1^{\alpha_1} \ldots x_n^{\alpha_n}}, 
\{ f=0 \} \right) 
&\longmapsto& \left\{ \frac{\partial f}{\partial 
x_1^{\alpha_1} \ldots x_n^{\alpha_n}} = 0 \right\}
\end{eqnarray*}
where $\alpha_1 + \ldots + \alpha_n = d$.
By the Hopf Theorem (see for instance \cite{ACGH:85}, p.~108), we have 
$$
h_{m+d}(p) \ge h_m(p) + {{n-1+d} \choose {n-1}} - 1. 
$$
\end{Remark} 

\begin{Theorem}\label{main}
Let $X$ be a smooth projective variety of dimension $n$. Assume that $X$ 
satisfies  
$$
N \ge (m-1)\left[{{n+m-1}\choose{n-1}}-n-1 \right] - \sum_{d=1}^{m-2} 
{{n-1+d}\choose{n-1}}
$$
independent differential equations of order $m$ at a general point $p \in X$. 
Then either $X \subset \PP^{{{n+m}\choose{n}}-N-1}$ or $X$ is covered by 
infinitely many subvarieties $Y$ of dimension at least $n-k \ge 1$ such 
that $Y \subset \PP^{{{n+m}\choose{n}}-N-k-m}$.
\end{Theorem}

\proof We split the proof into three cases, according to the possible values 
of $h_m(p)$. 

If $h_m(p) = {{n+m}\choose{n}} - {{n+m-1}\choose{n}}
= {{n+m-1}\choose{n-1}}$, then $T(m,p,X) = T(m-1,p,X)$. From \cite{ChiCil:02}, 
Proposition~2.3, it follows that $X \subset \PP^{{{n+m}\choose{n}}-N-1}$. 

If $h_m(p) = {{n+m-1}\choose{n-1}} - k$ with 
$1 \le k \le n-1$, then $\dim T(m,p,X) = \dim T(m-1,p,X) + k$. From 
Theorem~\ref{bompiani} it follows that either $X \subset \PP^{{{n+m}\choose{n}}
-N-1}$ or $X$ is covered by infinitely many subvarieties $Y$ of dimension 
at least $n-k$ such that $Y \subset \PP^{{{n+m}\choose{n}}-N-k-m}$. 

Assume finally $h_m(p) \le {{n+m-1}\choose{n-1}} - n$. 
From Remark~\ref{hopf} it follows that $h_{m-d} \le {{n+m-1}\choose{n-1}} - n 
- {{n-1+d}\choose{n-1}} + 1$ for every $d \le m-2$. Therefore   
$X$ satisfies at most   
$(m-1){{n+m-1}\choose{n-1}} -(m-1)n + (m-2) - \sum_{d=1}^{m-2} 
{{n-1+d}\choose{n-1}}$ equations and this contradiction ends the proof.

\qed

\section{}\label{s2}
Let $X$ be an integral $n$-dimensional projective variety defined over an 
algebraically closed field $\K$ of arbitrary characteristic.  
Fix a proper closed subset $T$ of $X$, $L\in \mbox{Pic}(X)$, and 
$V \subseteq H^0(X,L)$ such that $V$ spans $L$. Hence $V$ induces a morphism 
$\phi _V: X \to \PP(V^\ast ) \cong \PP^r$, with $r:= 
\dim (V)-1$.
Set $U:= X\backslash T$. We always assume that $U$ is smooth and that 
$\phi _V\vert U$ is an embedding. We aim to find conditions
on $X$, $L$, $V$ and $T$ which force the existence of $q\in \phi _V(U)$ 
which is a hyperosculating point of the variety $\phi _V(X)$.
Of course if $T=\emptyset$, then $X$ is smooth and $\phi _V$ 
is an embedding. For any integer $m \ge 0$, 
let $P^m(L)$ be the sheaf of principal parts of order ar most $m$ of $L$ 
(\cite{p}, \S 2 and \S 6). Since $U$ is smooth and $n$-dimensional, 
$P^m(L)\vert U$ is locally free of rank ${{n+m} \choose {n}}$. 
There is a Taylor series map $a^m: V\otimes \mathcal {O}_X \to P^m(L)$ 
and the sheaves $\mbox{Im}(a^m)$ and $\mbox{Coker}(a^m)$ measure the 
osculating behaviour up to order $m$ of the variety $\phi _V(X) \subset 
\PP^r$. In particular, if $T(m,p,X)$ is the $m$-osculating space as defined 
in the previous section, we have $T(m,p,X)= \PP(\mbox{Im}(a^m(p))) 
\subseteq \PP(V^\ast)$. We say that $q\in U$ is a hyperosculating point of 
order at most $m$ for $V$ if the sheaf $\mbox{Coker}(a^m)$ is not locally free
at $q$; we use the convention that $q\in U$ is not a hyperosculating point of 
order at most $m$ for $V$ if $\mbox{Coker}(a^m)$ is the zero-sheaf in a 
neighborhood of $q$. Since $U$ is reduced, the set of hyperosculating points 
is a proper closed subset of $U$. 

\begin{Remark}
Let $X$ be a reduced algebraic variety and $\mathcal{F}$ a coherent sheaf on 
$X$. For every $ p \in X$ set $\alpha(p) = \dim_\K \mathcal{F}_p / 
\mathfrak{m}_p\mathcal{F}_p$, where $\mathcal{F}_p$ is the stalk of 
$\mathcal{F}$ at $p$ and $\mathfrak{m}_p$ is the maximal ideal of 
$\mathcal{O}_{X,p}$. It is easy to check that the function $\alpha: 
X \to \N$ is locally constant if and only if $\mathcal{F}$ is locally 
free. It follows that if $X \subset \PP^n$ and $q \in X_\reg$, then $q$ is 
hyperosculating if and only if there is $m > 0$ such that 
$\dim T(m,q,X) < \dim T(m,p,X)$ for a general $p \in X$.  
\end{Remark}

Set $x(m):= \mbox{rank}(\mbox{Im}(a_m))$. The main result of \cite{fkpt} 
can be formulated as follows:

\begin{Theorem}\label{projective}\emph{(\cite{fkpt})}
Assume $T=\emptyset$, $V = H^0(X,L)$ and $r+1 = {{n+m} \choose {m}}$. 
If $x(m)=r+1$ and there is no hyperosculating point of order at
most $m$, then $X \cong \PP^n$ and $L \cong 
\mathcal {O}_{\PP^n}(m)$.
\end{Theorem}

Our set-up suggests the following natural generalization:

\begin{Proposition}\label{a3}
Assume $X$ smooth, $r+1 = {{n+m}\choose{m}}= x(m)$, $\codim(T) \ge 2$ 
and that there is no hyperosculating point of order at most $m$ on $U$, 
then $X \cong \PP^n$ and $L \cong \mathcal {O}_{\PP^n}(m)$.
\end{Proposition}

\proof By assumption there is an everywhere injective map (with
locally free or zero cokernel) $a^m\vert U: V\otimes \mathcal {O}_U \to P^m(L)
\vert U$. Since $X$ is reduced, this implies the injectivity of $a^m$ as a map 
of sheaves. Since $X$ is smooth, $P^m(L)$ is locally free. An injective map
of sheaves between two locally free sheaves with the same rank is either an 
isomorphism or its cokernel is supported exactly by an effective Cartier 
divisor (use the determinant). Hence $a^m$ is an isomorphism. 
By Theorem~\ref{projective}, $X \cong \PP^n$ and $\phi _V$ is a Veronese 
embedding.

\qed

Next we introduce a deeper result which points to the same direction. 
For any vector bundle $E$ on $X$ (or on $U$), let $c_\ast (E)$ denote 
its total Chern class in the Chow ring of $X$ (resp. $U$).
Since $U$ is smooth, for every integer $t>0$ we have the following exact 
sequence on $U$:
\begin{equation}\label{eqb1}
0 \to S^t(\Omega ^1_U) \to P^t(L)\vert U \to P^{t-1}(L)\vert U \to 0
\end{equation}
Hence by induction on $m$ one can compute the Chern classes on $U$ 
(not on $X$) of $P^m(L)\vert U$ in terms of $c_1(L\vert U)$
and the Chern classes of $\Omega ^1_U$, i.e. the Chern classes of $U$.

\begin{Theorem}\label{a1}
Assume $x(m) = r+1 \le {{n+m}\choose{n}}\le r+n$, and that there is no 
hyperosculating point of order at most $m$ on $U$. 
Then $c_t(P^m(L)\vert U)=0$ for every integer $t$ such that 
${{n+m}\choose{n}}-r-1 < t \le n$.
\end{Theorem}

\proof Set $A_U:= \mbox{Coker}(a^m\vert U)$. Since $x(m)=r+1$, $a^m$ 
is generically injective. Since $X$ is reduced
and $P^m(L)\vert U$ has no torsion, it follows that $a^m$ is injective as a 
map of sheaves. Hence $\mbox{Im}(a^m\vert U) \cong V\otimes \mathcal {O}_U$
and we have an exact sequence on $U$:
\begin{equation}\label{eqb2}
0 \to  V\otimes \mathcal {O}_U \to P^m(L)\vert U \to A_U \to 0
\end{equation}
By assumption, $A_U$ is locally free with rank ${{n+m}\choose{n}}-r-1<n$. 
Hence $c_t(A_U) = 0$ for every integer $t$ such that
$n \ge t > \mbox{rank}(A_U)$. By (\ref{eqb2}), we have $c_i(P^m(L)\vert U) = 
c_i(A_U)$ for every $i$, so the proof is over.

\qed

We have also the following result.

\begin{Proposition}\label{a4}
Set $R:= \omega _U^{\otimes {{n+m}\choose{n+1}}}\otimes 
L^{\otimes {{n+m}\choose{n}}}\vert U$. Assume $r = {{n+m} \choose{n}}$ and that
$a^m$ is surjective at each point of $U$. Then $c_\ast (P^m(L)\vert U) = 
1/c_\ast (R)$ in the Chow group of $U$.
\end{Proposition}

\proof We have $R \cong \mbox{det}(P^m(L)\vert U)$ 
(use $m$ times (\ref{eqb1})).
Since $\mbox{dim}(V) = \mbox{rank}(P^m(L)\vert U)+1$ and $a^m\vert U$ is 
surjective, $\mbox{Ker}(a^m\vert U)$ is a line bundle on $U$. 
This line bundle is $R^\ast$ and we have the following exact sequence on $U$:
\begin{equation}\label{eqb3}
0 \to R^\ast \to V\otimes \mathcal {O}_U \to P^m(L)\vert U \to 0
\end{equation}
which gives the claimed relation between the Chern polynomials.

\qed

\section{}\label{s3}
Here we maintain the previous notation, but from now on we have to require 
$\mbox{char}(\K)=0$. First of all, we show that very mild assumptions force 
the existence of hyperosculating points. 

\begin{Theorem}\label{a5}
Let $x(m)=r+1 = {{n+m}\choose{m}}-1$ and set $J:= \omega _X^{\otimes m}\otimes
L^{\otimes (n+1)}$. If for some integer $y$ with $0 \le y \le n-2$ and 
a big and nef divisor $H$ on $X$ we have $J^{n-y}\cdot H^y > 0$ (intersection 
product of $n$ divisors), then $\phi _V(X)$ must have hyperosculating points
of order at most $m$.
\end{Theorem}

\proof We divide the proof into three steps.

\quad {\emph {Step 1)}} Here we assume the existence of an integral curve 
$D \subset X$ such that $-n-1 \le \omega _X\cdot D
< 0$ and $D\cdot L < m$. Hence the linear span $\langle \phi _V(D)\rangle$ 
of $\phi _V(D)$ has dimension at most $m-1$. Chosen
$p\in \phi _V(D)_{reg}$ and a local coordinate corresponding to the 
direction of $\phi _V(D)$, we recover a linear relation about the first
$m$ partial derivatives in that direction.

\quad {\emph {Step 2)}} 
Set $R:= \omega _X^{\otimes {{n+m}\choose{n+1}}}\otimes L^{\otimes 
{{n+m}\choose{n}}}$ and notice that a positive multiple of $R$ is 
isomorphic to a positive multiple of $J$.
By Step 1) we may assume that for every integral curve $D \subset X$ such 
that $-n-1 \le \omega _X\cdot D < 0$
we have $R\cdot D \ge 0$. For every integral curve $E\subset X$ such that 
$\omega _X\cdot E \ge 0$ we have $R\cdot E>0$ because
$L$ is ample. Set $\varepsilon := m/(n+1)$ and 
let $N_\varepsilon (X,L) := \{Z\in N(X)_\R: \omega_X \cdot Z +  \varepsilon 
L \cdot Z \ge 0\}$, where $N(X)_\R$ is the tensor product with $\R$ of the 
free group of $1$-cycles on $X$ modulo numerical equivalence (\cite{m}).
By the very definition of $\varepsilon$, we have $R \cdot Z \ge 0$ for every 
$Z \in N_\varepsilon (X,L)$. By the Cone Theorem
(\cite{m}, Theorem 1.4), $R$ is nef. Since $R$ is a rational 
multiple of $J$, we have $R^{n-y}\cdot H^y > 0$. From \cite{bs},
Lemma 2.2.7, it follows that $h^1(X,R^\ast )=0$.

\quad {\emph {Step 3)}} 
Assume by contradiction that there are no hyperosculating 
points of order at most $m$. 
As in the proof of Proposition \ref{a4}, we have
an exact sequence
\begin{equation}\label{eqb3.1}
0 \to \mathcal {O}_X^{\oplus x} \to P^m(L) \to R \to 0
\end{equation}
where $x:={{n+m}\choose{m}}-1$.
By Step 2), we have $h^1(X,R^\ast )=0$. 
Hence the exact sequence (\ref{eqb3.1}) splits and we obtain
$P^m(L) \cong \mathcal {O}^{\oplus x}\oplus R$ with $x:= {{n+m}\choose{n}} 
-1$. Let $C$ be the intersection of $n-1$ general
members of $\vert L\vert$. Hence $(L^\ast )^{\oplus (n-1)}$ is the 
conormal module of $C$ in $X$ and $(\Omega ^1_X\vert C)\otimes L)$
is an extension of $\omega _C\otimes L$ by $\mathcal {O}_C^{\oplus (n-1)}$. 
Hence $S^m(\Omega^1_X )\otimes L^{\otimes m}\vert C$
is an extension of an ample rank ${{n+m-1}\choose{m}} - {{n+m-2}\choose{m}}$ 
vector bundle $F$ by the trivial rank ${{n+m-2}\choose{m}}$
vector bundle. By (\ref{eqb1}), there is an injection
$j: S^m(\Omega^1_X )\otimes L^{\otimes m}\vert C \to P^m(L)\vert C \cong
\mathcal {O}_C^{\oplus x}\oplus (R\vert C)$. At least ${{n+m-2}\choose{m}}-1$ 
of the trivial factors of $S^m(\Omega^1_X)\otimes L^{\otimes m}\vert C$
are mapped isomorphically onto some of the trivial factors of $P^m(L)\vert C$. 
Hence $j$ induces a map $u: F \to \mathcal {O}_C^{
\oplus x}$ with $\mbox{rank}(u) \ge \mbox{rank}(F)-1>0$, 
in contradiction with the ampleness of $F$.

\qed

\begin{Remark}\label{a6}
Let $C \subset \PP^r$ be an integral non-degenerate curve and 
$f: X \to C$ its normalization
map. Take $L:= f^\ast (\mathcal {O}_C(1))$ and 
$V:= f^\ast (H^0(\PP^r,\mathcal {O}_{{\bf {P}}^r}(1)))$.
By the Brill - Segre formula (see e.g. \cite{l}, p. 54), there are no 
hyperosculating points if and only if $C$ is a rational normal curve.
\end{Remark}

Sometimes it is also possible to bound the drop of the dimension of 
the osculating space at a hyperosculating point. For instance, consider
the following generalization of a result by Lanteri (see \cite{lan}, 
Theorem~B): 

\begin{Proposition}\label{lanteri}
Assume that $X$ is a linear $\PP^{n-1}$-bundle over a smooth curve. Then 
$\dim T(2,p,X) \ge n+1$ for every $p \in X$. 
\end{Proposition}

\proof We follow an argument provided by Lanteri in \cite{lan} for 
the $2$-dimensional case. Since $\dim V = \dim \Ker(\alpha^m(p)) + 
\dim \im(\alpha^m(p))$, we have 
$$
\dim T(m,p,X) + \dim \vert V - (m+1)p \vert = r - 1
$$ 
for every $m \ge 1$, and we deduce that 
$$
\dim T(2,p,X) = n + \codim(\vert V - 3p \vert, \vert V - 2p \vert).
$$
Therefore it is sufficient to show that $\vert V - 3p \vert \ne 
\vert V - 2p \vert$ for every $p \in X$. In order to do so, 
let $F(p)$ be the fiber of $X$ through $p$ and notice that 
$$
\vert V - 3p \vert = 2 F(p) + \vert V - 2 F(p) - p \vert,
$$
as it easily follows from the fact that $D^{n-1}.F(p)=1$ for every 
$D \in \vert V \vert$. 
If $\vert V - 2p \vert = 2 F(p) + \vert V - 2 F(p) - p \vert$
for some point $p \in X$, then every hyperplane tangent to $X$ 
at $p$ is tangent along the whole fiber $F(p)$. In particular,
the tangent space to $X$ is constant along a positive dimensional 
subvariety, in contradiction with Zak's Theorem on the finiteness 
of the Gauss map (see for instance \cite{zak}, Chapter~I, Corollary~2.8). 

\qed

\vspace{0.5cm}

\noindent
Edoardo Ballico \newline
Universit\`a degli Studi di Trento \newline
Dipartimento di Matematica \newline
Via Sommarive 14 \newline
38050 Povo (Trento) \newline
Italy \newline
e-mail: ballico@science.unitn.it

\vspace{0.5cm}

\noindent
Claudio Fontanari \newline
Universit\`a degli Studi di Trento \newline
Dipartimento di Matematica \newline
Via Sommarive 14 \newline
38050 Povo (Trento) \newline
Italy \newline
e-mail: fontanar@science.unitn.it

\end{document}